\documentclass[review]{elsarticle}
\usepackage{amsmath}
\usepackage{tikz}
\usepackage{verbatim}
\usepackage{amssymb}
\usepackage{mathtools}
\usepackage{subcaption}
\usepackage{pgfplots}
\usepackage{booktabs}
\usepackage{tcolorbox}
\usepackage{pgfplotstable}
\usepackage{lineno,hyperref}
\usepackage{xintexpr}

\pgfplotsset{compat=1.11,
    /pgfplots/ybar legend/.style={
    /pgfplots/legend image code/.code={%
       \draw[##1,/tikz/.cd,yshift=-0.25em]
        (0cm,0cm) rectangle (3pt,0.8em);},
   },
}

\modulolinenumbers[5]

\begin{document}

\begin{frontmatter}
\title{Combining $p$-multigrid and multigrid reduced in time methods to obtain a scalable solver for Isogeometric Analysis}
\tnotetext[mytitlenote]{Fully documented templates are available in the elsarticle package on \href{http://www.ctan.org/tex-archive/macros/latex/contrib/elsarticle}{CTAN}.}

\author[mymainaddress]{Roel Tielen\corref{mycorrespondingauthor}}
\cortext[mycorrespondingauthor]{Corresponding author}
\ead{r.p.w.m.tielen@tudelft.nl}

\author[mymainaddress]{Matthias M\"oller}
\author[mymainaddress]{Cornelis Vuik}

\address[mymainaddress]{Delft University of Technology, Mekelweg 4, 2628 CD, Delft}

\begin{abstract}
Isogeometric Analysis (IgA) \cite{Hughes2005} has become a viable alternative to the Finite Element Method (FEM) and is typically combined with a time integration scheme within the method of lines for time-dependent problems. However, due to a stagnation of processor's clock speeds, traditio\-nal (i.e. sequential) time integration schemes become more and more the bottleneck within these large-scale computations, which lead to the development of parallel-in-time methods like the Multigrid Reduced in Time (MGRIT) method \cite{Falgout2014}. 

Recently, MGRIT has been succesfully applied by the authors in the context of IgA showing convergence independent of the mesh width, approximation order of the B-spline basis functions and time step size for a variety of benchmark problems. However, a strong dependency of the CPU times on the approximation order was visible when a standard Conjugate Gradient method was adopted for the spatial solves within MGRIT. In this paper we combine MGRIT with a state-of-the-art solver (i.e. a $p$-multigrid method \cite{Tielen2019}), specifically designed for IgA, thereby significantly reducing the overall computational costs of MGRIT. Furthermore, we investigate the performance of MGRIT and its scalability on modern copmuter architectures. 
\end{abstract}


\end{frontmatter}


\section{Introduction}
Since its introduction in \cite{Hughes2005}, Isogeometric Analysis (IgA) has become more and more a viable alternative to the Finite Element Method (FEM). Within IgA, the same building blocks (i.e. B-splines and NURBS) as in Computer Aided Design (CAD) are adopted, which closes the gap between CAD and FEM. In particular, the use of high-order splines results in a highly accurate represention of (curved) geometries and has shown to be advantageous in many applications, like structural mechanics \cite{Cottrell2006,Kiendl2009,Kiendl2010,Benson2010}, solid and fluid dynamics \cite{Bazilevs2006,Moutsanidis2020,Gan2017,Tielen2017a} and shape optimization \cite{Wall2008,Qian2010,Seo2010,Li2011}. Finally, the accuracy per degree of freedom (DOF) compared to FEM is significantly higher with IgA \cite{Hughes2008}.

For time-dependent partial differential equations (PDEs), IgA is typically combined with a traditional time integration scheme within the method of lines. Here, the spatial variables are discretizated by adopting IgA, after which the resulting system of ordinary differential equations (ODEs) is integrated in time. However, as with all traditional time integration schemes, the latter part becomes more and more the bottleneck in numerical simulations. When the spatial resolution is increased to improve accuracy, a smaller time step size has to be chosen to ensure stability of the overall method. As clock speeds are no longer increasing, but the core count goes up, the parallelizability of the entire calculation process becomes more and more important to obtain an overall efficient method. As traditional time integration schemes are sequential by nature, new parallel-in-time methods are needed to resolve this problem.   

The Multigrid Reduced in Time (MGRIT) method \cite{Falgout2014} is a parallel-in-time algorithm based on multigrid reduction (MGR) techniques \cite{Ries1983}. In contrast to space-time methods, in which time is considered as an extra spatial dimension, sequential time stepping is still necessary within MGRIT. Space-time methods have been combined in the literature with IgA \cite{Langer2016,Takizawa2017,Hofer2018,Hofer2019}. Although very successful, a drawback of such methods is the fact that they are more intrusive on existing codes, while MGRIT just requires a routine to integrate the fully discrete problem from one time instance to the next. Over the years, MGRIT has been studied in detail (see \cite{Bal2005,Gander2007,Gander2008,Dobrev2017,Southworth2019}) and applied to a variety of problems, including those arising in optimization \cite{Gnther2018b,Gnther2018} and power networks \cite{Lecouvez,Schroder}.        

Recently, the authors applied MGRIT in the context of IgA for the first time in the literature \cite{Tielen2021}. Here, MGRIT showed convergence for a variety of twodimensional benchmark problems independent of the mesh width $h$, the approximation order $p$ of the B-spline basis functions and the number of time steps $N_t$. However, as a standard Conjugate Gradient method was adopted for the spatial solves within MGRIT, a significant dependency of the CPU times on the approximation order was visible. Furthermore, the parallel performance of MGRIT was investigated for a limited number of cores.  

In this paper, we combine MGRIT with a state-of-the-art $p$-multigrid method \cite{Tielen2019} to solve the linear systems arising within MGRIT. CPU timings show that the use of such a solver significantly improves the overall performance of MGRIT, in particular for higher values of $p$. Furthermore, the parallel performance of the resulting MGRIT method (i.e. strong and weak scalability) is investigated on modern computer architectures.  

This paper is structured as follows: In Section \ref{sec:model}, a two-dimensional model problem and its spatial and temporal discretization are considered. The MGRIT algorithm is then described in Section \ref{sec:method}. In Section \ref{sec:pmg}, the adopted $p$-multigrid method and its components are presented in more detail. In Section \ref{sec:results}, numerical results obtained for the considered model problem are analyzed for different values of the mesh width, approximation order and the number of time steps. Furthermore, weak and strong scaling studies of MGRIT when combined with IgA are performed. Finally, conclusions are drawn in Section \ref{sec:conclusions}.  

\section{Model problem and discretization}
\label{sec:model}

As a model problem, we consider the transient diffusion equation:
\begin{eqnarray} \label{eq:heat}
\partial_t u(\mathbf{x},t) - \kappa \Delta u(\mathbf{x},t) = f(\mathbf{x}), \quad \mathbf{x} \in \Omega, t \in [0,T].
\end{eqnarray}

Here, $\kappa$ denotes a constant diffusion coefficient, $\Omega \subset \mathbb{R}^d$ a simply connected, Lipschitz domain in $d$ dimensions and $f \in L^2(\Omega)$ a source term. The above equation is complemented by initial conditions and homogeneous Dirichlet boundary conditions:
\begin{eqnarray}
u(\mathbf{x},0) &=& u^0(\mathbf{x}), \quad \mathbf{x} \in \Omega, \\
u(\mathbf{x},t) &=& 0, \quad \quad \  \ \mathbf{x} \in \partial \Omega, t \in [0,T].
\end{eqnarray}

First, we discretize Equation \eqref{eq:heat} by dividing the time interval $[0,T]$ in $N_t$ subintervals of size $\Delta t$ and applying the $\theta$-scheme to the temporal derivative, which leads to the following equation to be solved at every time step:
\begin{eqnarray} 
\frac{u(\mathbf{x})^{k+1}-u(\mathbf{x})^{k}}{\Delta t} = \kappa \theta \Delta u(\mathbf{x})^{k+1} + \kappa (1-\theta) u(\mathbf{x})^k + f(\mathbf{x}), 
\end{eqnarray}
for $\mathbf{x} \in \Omega$ and  $k = 0,\ldots,N_t$. Depending on the choice of $\theta$, this scheme leads to the backward Euler ($\theta=1$), forward Euler ($\theta=0$) or Crank-Nicolson ($\theta=0.5$) method, which will all be adopted throughout this paper. By rearranging the terms, the discretized equation can be written as follows: 
\begin{eqnarray} \label{eq:semidiscrete}
u(\mathbf{x})^{k+1} - \kappa \Delta t\theta \Delta u(\mathbf{x})^{k+1} = u(\mathbf{x}^k) + \kappa \Delta t(1-\theta) u(\mathbf{x})^k + \Delta t f(\mathbf{x}).
\end{eqnarray}

To obtain the variational formulation, let $\mathcal{V} = H^1_0(\Omega)$ be the space of functions in the Sobolev space  $H^1(\Omega)$ that vanish on the boundary $\partial \Omega$. Equation~\eqref{eq:semidiscrete} is multiplied with a test function $v \in \mathcal{V}$ and the result is then integrated over the domain $\Omega$:
\begin{eqnarray} 
\int_{\Omega} u^{k+1}  v  - \kappa \Delta t \theta \Delta u^{k+1} v \text{d}\Omega  = \int_{\Omega} u^{k}  v  + \kappa \Delta t (1-\theta)  \Delta u^{k} v +  \Delta t f v \ \text{d}\Omega. 
\end{eqnarray}

Applying integration by parts on the second term on both sides of the equation results in 
\begin{eqnarray}  \label{eq:weakform}
\int_{\Omega} u^{k+1} v + \kappa \Delta t \theta \nabla u^{k+1} \cdot \nabla v \ \text{d}\Omega  = \int_{\Omega} u^{k+1} v - \kappa \Delta t (1-\theta) \nabla u^k \cdot \nabla v + \Delta tf v \ \text{d}\Omega, 
\end{eqnarray}
where the boundary integral integral vanishes since $v = 0$ on $\partial \Omega$. To parameterize the physical domain $\Omega$, a geometry function $\mathbf{F}$ is then defined, describing an invertible mapping to connect the parameter domain $\Omega_0 = (0,1)^d$ with the physical domain $\Omega$:
\begin{equation}
    \mathbf{F} = \Omega_0 \rightarrow \Omega, \quad \mathbf{F}(\boldsymbol{\xi}) = (\mathbf{x}).
\end{equation}
 
Provided that the physical domain $\Omega$ is topologically equivalent to the unit square, the geometry can be described by a single geometry function $\mathbf{F}$. In case of more complex geometries, a family of functions $\mathbf{F}^{(m)}$ ($m = 1,\ldots, K$) is defined and we refer to $\Omega$ as a multipatch geometry consisting of $K$ patches. For a more detailed description of the spatial discretization in IgA and multipatch constructions, the authors refer to chapter $2$ of \cite{Hughes2005}.

At each time step, we express $u$ in Equation \eqref{eq:weakform} by a linear combination of multivariate B-spline basis functions of order $p$. Multivariate B-spline basis functions are defined as the tensor product of univariate B-spline basis functions $\phi_{i,p}$ $(i=1,\ldots,N)$ which are uniquely defined on the parameter domain $(0,1)$ by an underlying knot vector $\Xi = \{ \xi_1, \xi_2, \ldots , \xi_{N+p}, \xi_{N+p+1} \}$. Here, $N$ denotes the number of B-spline basis functions and $p$ the approximation order. Based on this knot vector, the basis functions are defined recursively by the Cox-de Boor formula \cite{Boor1978}, starting from the constant ones
\begin{eqnarray}
\phi_{i,0}(\xi) = \begin{cases} 1   \hspace{2.25 cm} \text{if} \hspace{0.2cm}\xi_i \leq \xi < \xi_{i+1}, \\  0 \hspace{2.25cm} \text{otherwise.}  \end{cases}
\end{eqnarray}

Higher-order B-spline basis functions of order $p>0$ are then defined recursively

\begin{eqnarray}
\phi_{i,p}(\xi) = \frac{\xi - \xi_i}{\xi_{i+p}-\xi_i} \phi_{i,p-1}(\xi) +  \frac{\xi_{i+p+1} - \xi}{\xi_{i+p+1}-\xi_{i+1}} \phi_{i+1,p-1}(\xi).
\end{eqnarray}

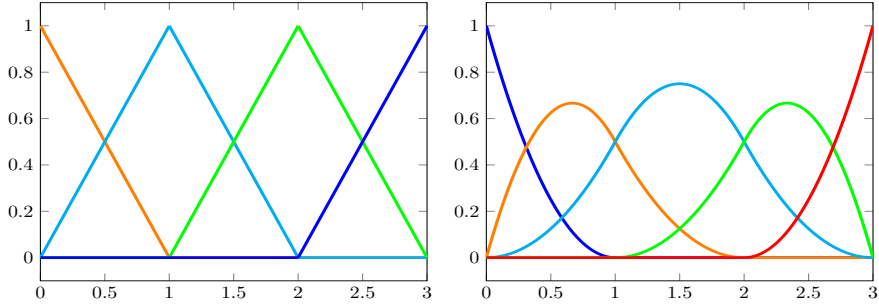
\begin{figure}[ht!]
\centering
\begin{tikzpicture}[xscale=0.75,yscale=0.65]
\begin{axis}[xmin=0,xmax=3,ymin=-0.1,ymax=1.1, samples=50]
  \addplot[orange, ultra thick, domain=0:1] {1-x};
  \addplot[orange, ultra thick, domain=1:3] {0};
  \addplot[cyan, ultra thick, domain=0:1] {x};
  \addplot[cyan, ultra thick, domain=1:2] {2-x};
  \addplot[cyan, ultra thick, domain=2:3] {0};
  \addplot[green, ultra thick, domain=0:1] {0};
  \addplot[green, ultra thick, domain=1:2] {x-1};
  \addplot[green, ultra thick, domain=2:3] {3-x};
  \addplot[blue, ultra thick, domain=0:2] {0};
  \addplot[blue, ultra thick, domain=2:3] {x-2};
\end{axis}
\end{tikzpicture}
\begin{tikzpicture}[xscale=0.75,yscale=0.65]
\begin{axis}[xmin=0,xmax=3,ymin=-0.1,ymax=1.1, samples=50]
  \addplot[blue, ultra thick, domain=0:1] {(1-x)*(1-x)};
  \addplot[blue, ultra thick, domain=1:3] {0};
  \addplot[orange, ultra thick, domain=0:1] {x*(1-x)+0.5*(2-x)*x};
  \addplot[orange, ultra thick, domain=1:2] {0.5*(2-x)*(2-x)};
  \addplot[orange, ultra thick, domain=2:3] {0};
  \addplot[cyan, ultra thick, domain=0:1] {0.5*x*x+0.5*(3-x)*0};
  \addplot[cyan, ultra thick, domain=1:2] {0.5*x*(2-x)+0.5*(3-x)*(x-1)};
  \addplot[cyan, ultra thick, domain=2:3] {0.5*x*0+0.5*(3-x)*(3-x)};
  \addplot[green, ultra thick, domain=0:1] {0};
  \addplot[green, ultra thick, domain=1:2] {0.5*(x-1)*(x-1)};
  \addplot[green, ultra thick, domain=2:3] {0.5*(x-1)*(3-x)+(3-x)*(x-2)};
  \addplot[red, ultra thick, domain=0:2] {0};
  \addplot[red, ultra thick, domain=2:3] {(x-2)*(x-2)};
\end{axis}
\end{tikzpicture}
\caption{Linear and quadratic B-spline basis functions based on the knot vectors $\Xi_1 = \{ 0, 0, 1, 2, 3,3 \}$ and $\Xi_2 = \{ 0, 0, 0, 1, 2, 3,3,3 \}$, respectively.}
\label{fig:bspline}
\end{figure}

The resulting B-spline basis functions $\phi_{i,p}$ are non-zero on the interval $[\xi_i,\xi_{i+p+1})$ and possess the partition of unity property. Furthermore, the basis functions are $C^{p-m_i}$-continuous, where $m_i$ denotes the multiplicity of knot $\xi_i$. Throughout this paper, we consider a uniform knot vector with knot span size $h$, where the first and last knot are repeated $p+1$ times. As a consequence, the resulting B-spline basis functions are $C^{p-1}$ continuous and interpolatory at both end points. Figure \ref{fig:bspline} illustrates both linear and quadratic B-spline basis functions based on such a knot vector. 

As mentioned previously, the tensor product of univariate B-spline basis functions is adopted for the multi-dimensional case. Denoting the total number of multivariate B-spline basis functions $\Phi_{i,p}$ by $N_{\rm dof}$, the solution $u$ is thus approximated as follows:
\begin{eqnarray}
u(\mathbf{x}) \approx u_{h,p}(\mathbf{x}) = \sum_{i=1}^{N_{\rm dof}} u_i(t) \Phi_{i,p}(\mathbf{x}), \quad u_{h,p} \in \mathcal{V}_{h,p}.
\end{eqnarray} 

\noindent Here, the spline space $\mathcal{V}_{h,p}$ is defined, using the inverse of the geometry mapping $\mathbf{F}^{-1}$ as pull-back operator, as follows:
\begin{eqnarray}
\mathcal{V}_{h,p} := \text{span}\left  \{ \Phi_{i,p} \circ \mathbf{F}^{-1} \right \}_{i=1, \ldots,N_{\rm dof}}.
\end{eqnarray}

\noindent By setting $v = \Phi_{j,p}$, Equation \eqref{eq:weakform} can be written as follows:
\begin{eqnarray} \label{eq:los}
\left ( \mathbf{M} + \kappa \Delta t \theta \mathbf{K}  \right ) \mathbf{u}^{k+1}= \left ( \mathbf{M} - \kappa \Delta t (1-\theta) \mathbf{K} \right ) \mathbf{u}^{k} + \Delta t \mathbf{f}, \ k = 0,\ldots,N_t, 
\end{eqnarray}

\noindent where $\mathbf{M}$ and $\mathbf{K}$ denote the mass and stiffness matrix, respectively:
\begin{equation}
\mathbf{M}_{i,j} = \int_{\Omega}  \Phi_{i,p} \Phi_{j,p} \ \text{d}\Omega, \qquad \mathbf{K}_{i,j} = \int_{\Omega} \nabla \Phi_{i,p} \cdot \nabla \Phi_{j,p} \ \text{d}\Omega. 
\end{equation}

\section{Multigrid Reduced in Time}
\label{sec:method}

A traditional (i.e. sequential) time integration scheme would solve Equation \eqref{eq:los} for $k=0,\ldots,N_t$ to obtain the numerical solution at each time instance. In this paper, however, we apply MGRIT to solve Equation \eqref{eq:los} parallel-in-time. For the ease of notation, we set $\theta=1$ throughout the remainder of this section. Let $\Psi = \left ( \mathbf{M} + \kappa \Delta t \mathbf{K} \right )^{-1}$ denote the inverse of the left hand side operator.  Then, Equation \eqref{eq:los} can be written as follows:
\begin{eqnarray}
\mathbf{u}^{k+1}&=& \Psi \left ( \mathbf{M} \mathbf{u}^{k} + \Delta t \mathbf{f} \right ), \ k = 0,\ldots,N_t,  \\
                &=& \Psi \mathbf{M} \mathbf{u}^{k} + \mathbf{g}^{k+1}, \ k = 0,\ldots,N_t, 
\end{eqnarray} 
where $\mathbf{g}^{k+1} = \Psi \Delta t \mathbf{f}$. Setting $\mathbf{g}^0$ equal to the initial condition $u^0(\mathbf{x})$ projected on the spline space $\mathcal{V}_{h,p}$, the time integration method can be written as a linear system of equations:
\begin{eqnarray} \label{eq:mgrit}
\mathbf{A} \mathbf{u}= \begin{bmatrix} I &  &  & \\ -\Psi \mathbf{M} & I & &  \\   & \ddots &  \ddots &  \\   &  & -\Psi \mathbf{M} & I \end{bmatrix} \begin{bmatrix} \mathbf{u}^0 \\ \mathbf{u}^1 \\ \vdots \\ \mathbf{u}^{N_t} \end{bmatrix} = \begin{bmatrix} \mathbf{g}^0 \\ \mathbf{g}^1 \\ \vdots \\ \mathbf{g}^{N_t} \end{bmatrix} = \mathbf{g}.
\end{eqnarray}

A sequential time integration scheme would correspond to a block-forward solve of this linear system of equations. Here, we first introduce the two-level MGRIT method, showing similarities with the well-known parareal algorithm \cite{lions2001}. In fact, it can be shown that both methods are equivalent, assuming a specific choice of relaxation \cite{Gander2007}. Then, the multilevel variant of MGRIT will be presented in more detail.

\subsection*{Two-level MGRIT method}

\begin{figure}[h!]
\begin{center}
\begin{tikzpicture}
\draw (0,0) -- (10,0);
\draw(0,-0.5) -- (0,0.5);
\draw(0.5,-0.25) -- (0.5,0.25);
\draw(1.0,-0.25) -- (1.0,0.25);
\draw(1.5,-0.25) -- (1.5,0.25);
\draw(2.0,-0.25) -- (2.0,0.25);
\draw(2.5,-0.5) -- (2.5,0.5);
\draw(3.0,-0.25) -- (3.0,0.25);
\draw(3.5,-0.25) -- (3.5,0.25);
\draw(4.0,-0.25) -- (4.0,0.25);
\draw(4.5,-0.25) -- (4.5,0.25);
\draw(5,-0.5) -- (5,0.5);
\draw(5.5,-0.25) -- (5.5,0.25);
\draw(6.0,-0.25) -- (6.0,0.25);
\draw(6.5,-0.25) -- (6.5,0.25);
\draw(7.0,-0.25) -- (7.0,0.25);
\draw(7.5,-0.5) -- (7.5,0.5);
\draw(8.0,-0.25) -- (8.0,0.25);
\draw(8.5,-0.25) -- (8.5,0.25);
\draw(9.0,-0.25) -- (9.0,0.25);
\draw(9.5,-0.25) -- (9.5,0.25);
\draw(10,-0.5) -- (10,0.5);

\draw (0,0.7) node{$T_0$};
\draw (2.5,0.7) node{$T_1$};
\draw (3.75,0.7) node{$\cdots$};
\draw (10,0.7) node{$T_{N_t/m}$};
\draw (0,-0.8) node{$t_0$};
\draw (0.5,-0.8) node{$t_1$};
\draw (1.5,-0.8) node{$\cdots$};
\draw (2.5,-0.8) node{$t_m$};
\draw (10.0,-0.8) node{$t_{N_t}$};
\draw [thick, black,decorate,decoration={brace,amplitude=5pt,mirror},xshift=0.4pt,yshift=-0.4pt](8.0,-0.3) -- (8.5,-0.3) node[black,midway,yshift=-0.5cm] {\footnotesize $\Delta t_F$};
\draw [thick, black,decorate,decoration={brace,amplitude=5pt},xshift=0.4pt,yshift=-0.2pt](5.0,0.6) -- (7.5,0.6) node[black,midway,yshift=0.5cm] {\small $\Delta t_C = m \Delta t_F$};
\end{tikzpicture}
\end{center} 
\caption{Coarse and fine temporal mesh from $0$ to $T$.}
\label{fig:parareal}
\end{figure}
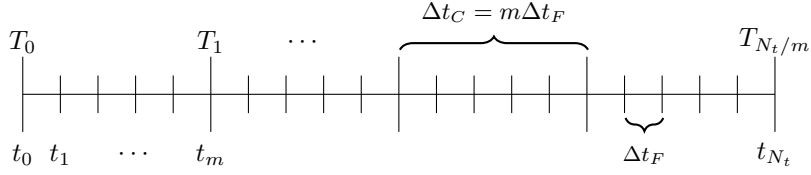

The two-level MGRIT method combines the use of a cheap coarse-level time integration method with an accurate more expensive fine-level one which can be performed in parallel. That is, the linear system of equations given by Equation \eqref{eq:mgrit} can be solved iteratively by introducing a coarse temporal mesh with time step size $\Delta t_C = m \Delta t_F$. Here, $\Delta t_F$ coincides with the $\Delta t$ from the previous sections and $m$ denotes the coarsening factor. Figure \ref{fig:parareal} illustrates both the fine and coarse temporal discretization.


It can be observed that the solution of Equation \eqref{eq:mgrit} at times $T_0, T_1, \ldots, T_{N_t/m}$ satisfies the following system of equations:
\begin{eqnarray} \label{eq:mgrit_coarse}
\mathbf{A}_{\Delta} \mathbf{u}_{\Delta}= \begin{bmatrix} I &  &  & \\ -(\Psi \mathbf{M})^m & I & &  \\   & \ddots &  \ddots &  \\   &  & -(\Psi \mathbf{M})^m & I \end{bmatrix} \begin{bmatrix} \mathbf{u}^0_{\Delta} \\ \mathbf{u}^1_{\Delta} \\ \vdots \\ \mathbf{u}^{N_t/m}_{\Delta} \end{bmatrix} = \begin{bmatrix} \mathbf{g}^0_{\Delta} \\ \mathbf{g}^1_{\Delta} \\ \vdots \\ \mathbf{g}^{N_t/m}_{\Delta} \end{bmatrix} = \mathbf{g}_{\Delta} .
\end{eqnarray}

Here, $\mathbf{u}_{\Delta}^j = \mathbf{u}^{jm}$ and the vector $\mathbf{g}_{\Delta}$ is given by the original vector $\mathbf{g}$ multiplied by a restriction operator:

\begin{eqnarray}
\mathbf{g}_{\Delta} = \begin{bmatrix} I & & & & & & & & & \\  & (\Psi \mathbf{M})^{m-1} & \cdots & \Psi \mathbf{M}  & I & & & & & \\ & & & & & \ddots & & & & \\ &  & & &  & & (\Psi \mathbf{M})^{m-1} & \cdots & \Psi \mathbf{M}& I \end{bmatrix}  \begin{bmatrix} \mathbf{g}^0 \\ \mathbf{g}^1 \\ \vdots \\ \mathbf{g}^{N_t/m} \end{bmatrix}
\end{eqnarray}

\begin{figure}[h!]
\begin{center}
\begin{tikzpicture}[xscale=0.75,yscale=0.75]
  \begin{axis}[xlabel=time,xmin=0,xmax=10,ylabel={Solution}] 
    \addplot[no marks,domain=0:10] {-0.1*(x-5)^2 + 3}; 
    \addplot[mark=*,blue,domain=0:1] {-0.1*(x-5)^2 + 3}; 
    \addplot[mark=*,blue,domain=2:3] {-0.1*(x-5)^2 + 3}; 
    \addplot[mark=*,blue,domain=4:5] {-0.1*(x-5)^2 + 3}; 
    \addplot[mark=*,blue,domain=6:7] {-0.1*(x-5)^2 + 3}; 
    \addplot[mark=*,blue,domain=8:9] {-0.1*(x-5)^2 + 3}; 
    \addplot[blue, mark=square] coordinates {(0,0.5)  (2,2.3)};
    \addplot[blue, mark=square] coordinates {(2,2.1)  (4,3.1)};
    \addplot[blue, mark=square] coordinates {(4,2.9)  (6,3.1)};
    \addplot[blue, mark=square] coordinates {(6,2.9)  (8,2.3)};
    \addplot[blue, mark=square] coordinates {(8,2.1)  (10,0.7)};
    \addplot[only marks,mark=square*,mark size=3pt, mark options={,color = red}] coordinates {(0,0.5) (2,2.1) (4,2.9) (6,2.9) (8,2.1) (10,0.5)};
    \end{axis}
\end{tikzpicture}
\centering \caption{Initial guess of the two-level MGRIT method and the (parallel) fine solutions.}
\label{fig:parareal2}
\end{center} 
\end{figure}
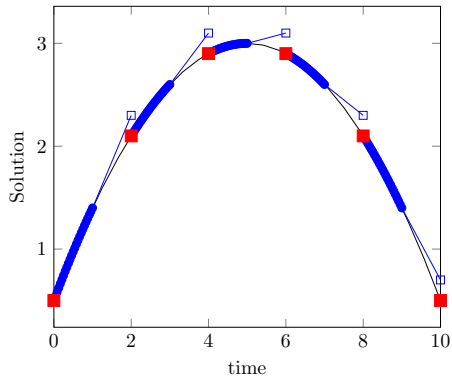

A two-level MGRIT method solves the coarse system given by Equation \eqref{eq:mgrit_coarse} iteratively and computes the fine values in parallel within each interval $(t_{jm},t_{(j+1)m-1})$. The coarse system is solved using the following residual correction scheme:
\begin{eqnarray}
\mathbf{u}^{(k+1)}_{\Delta} = \mathbf{u}^{(k)}_{\Delta} + \mathbf{B}^{-1}_{\Delta} \left ( \mathbf{g}_{\Delta} - \mathbf{A}_{\Delta} \mathbf{u}^{(k)}_{\Delta} \right ),
\end{eqnarray}

where $\mathbf{B}_{\Delta}$ is a coarse-level equivalent of the matrix $\mathbf{A}$. Here the fine values are computed in parallel, denoted by the action of operator $\mathbf{A}_{\Delta}$. This in contrast to the action of $\mathbf{B}_{\Delta}$ which typically is performed on a single processor. Figure \ref{fig:parareal2} illustrates how MGRIT computes the fine solution in parallel, based on an initial guess at the coarse time grid. 

\begin{figure}[h]
\begin{center}
\begin{tikzpicture}
\draw (0,0) -- (10,0);
\filldraw (0,0) circle (2.5pt) node[below,yshift=-0.2pt] {\footnotesize C};
\filldraw (0.5,0) circle (1pt) node[below,yshift=-0.2pt] {\footnotesize F};
\filldraw (1.0,0) circle (1pt) node[below,yshift=-0.2pt] {\footnotesize F};
\filldraw (1.5,0) circle (1pt) node[below,yshift=-0.2pt] {\footnotesize F};
\filldraw (2.0,0) circle (1pt) node[below,yshift=-0.2pt] {\footnotesize F};
\filldraw (2.5,0) circle (2.5pt) node[below,yshift=-0.2pt] {\footnotesize C};
\filldraw (3.0,0) circle (1pt) node[below,yshift=-0.2pt] {\footnotesize F};
\filldraw (3.5,0) circle (1pt) node[below,yshift=-0.2pt] {\footnotesize F};
\filldraw (4.0,0) circle (1pt) node[below,yshift=-0.2pt] {\footnotesize F};
\filldraw (4.5,0) circle (1pt) node[below,yshift=-0.2pt] {\footnotesize F};
\filldraw (5,0) circle (2.5pt) node[below,yshift=-0.2pt] {\footnotesize C};
\filldraw (5.5,0) circle (1pt) node[below,yshift=-0.2pt] {\footnotesize F};
\filldraw (6.0,0) circle (1pt) node[below,yshift=-0.2pt] {\footnotesize F};
\filldraw (6.5,0) circle (1pt) node[below,yshift=-0.2pt] {\footnotesize F};
\filldraw (7.0,0) circle (1pt) node[below,yshift=-0.2pt] {\footnotesize F};
\filldraw (7.5,0) circle (2.5pt) node[below,yshift=-0.2pt] {\footnotesize C};
\filldraw (8.0,0) circle (1pt) node[below,yshift=-0.2pt] {\footnotesize F};
\filldraw (8.5,0) circle (1pt) node[below,yshift=-0.2pt] {\footnotesize F};
\filldraw (9.0,0) circle (1pt) node[below,yshift=-0.2pt] {\footnotesize F};
\filldraw (9.5,0) circle (1pt) node[below,yshift=-0.2pt] {\footnotesize F};
\filldraw (10,0) circle (2.5pt) node[below,yshift=-0.2pt] {\footnotesize C};

\draw [->] (0.05,0.2) to [bend left=40] node[below] {} (0.45,0.2);
\draw [->] (0.55,0.2) to [bend left=40] node[below] {} (0.95,0.2);
\draw [->] (1.05,0.2) to [bend left=40] node[below] {} (1.45,0.2);
\draw [->] (1.55,0.2) to [bend left=40] node[below] {} (1.95,0.2);

\draw [->] (2.55,0.2) to [bend left=40] node[below] {} (2.95,0.2);
\draw [->] (3.05,0.2) to [bend left=40] node[below] {} (3.45,0.2);
\draw [->] (3.55,0.2) to [bend left=40] node[below] {} (3.95,0.2);
\draw [->] (4.05,0.2) to [bend left=40] node[below] {} (4.55,0.2);

\draw [->] (5.05,0.2) to [bend left=40] node[below] {} (5.45,0.2);
\draw [->] (5.55,0.2) to [bend left=40] node[below] {} (5.95,0.2);
\draw [->] (6.05,0.2) to [bend left=40] node[below] {} (6.45,0.2);
\draw [->] (6.55,0.2) to [bend left=40] node[below] {} (6.95,0.2);

\draw [->] (7.55,0.2) to [bend left=40] node[below] {} (7.95,0.2);
\draw [->] (8.05,0.2) to [bend left=40] node[below] {} (8.45,0.2);
\draw [->] (8.55,0.2) to [bend left=40] node[below] {} (8.95,0.2);
\draw [->] (9.05,0.2) to [bend left=40] node[below] {} (9.55,0.2);

\begin{scope}[yshift=-1.0cm]
\draw (0,0) -- (10,0);
\filldraw (0,0) circle (2.5pt) node[below,yshift=-0.2pt] {\footnotesize C};
\filldraw (0.5,0) circle (1pt) node[below,yshift=-0.2pt] {\footnotesize F};
\filldraw (1.0,0) circle (1pt) node[below,yshift=-0.2pt] {\footnotesize F};
\filldraw (1.5,0) circle (1pt) node[below,yshift=-0.2pt] {\footnotesize F};
\filldraw (2.0,0) circle (1pt) node[below,yshift=-0.2pt] {\footnotesize F};
\filldraw (2.5,0) circle (2.5pt) node[below,yshift=-0.2pt] {\footnotesize C};
\filldraw (3.0,0) circle (1pt) node[below,yshift=-0.2pt] {\footnotesize F};
\filldraw (3.5,0) circle (1pt) node[below,yshift=-0.2pt] {\footnotesize F};
\filldraw (4.0,0) circle (1pt) node[below,yshift=-0.2pt] {\footnotesize F};
\filldraw (4.5,0) circle (1pt) node[below,yshift=-0.2pt] {\footnotesize F};
\filldraw (5,0) circle (2.5pt) node[below,yshift=-0.2pt] {\footnotesize C};
\filldraw (5.5,0) circle (1pt) node[below,yshift=-0.2pt] {\footnotesize F};
\filldraw (6.0,0) circle (1pt) node[below,yshift=-0.2pt] {\footnotesize F};
\filldraw (6.5,0) circle (1pt) node[below,yshift=-0.2pt] {\footnotesize F};
\filldraw (7.0,0) circle (1pt) node[below,yshift=-0.2pt] {\footnotesize F};
\filldraw (7.5,0) circle (2.5pt) node[below,yshift=-0.2pt] {\footnotesize C};
\filldraw (8.0,0) circle (1pt) node[below,yshift=-0.2pt] {\footnotesize F};
\filldraw (8.5,0) circle (1pt) node[below,yshift=-0.2pt] {\footnotesize F};
\filldraw (9.0,0) circle (1pt) node[below,yshift=-0.2pt] {\footnotesize F};
\filldraw (9.5,0) circle (1pt) node[below,yshift=-0.2pt] {\footnotesize F};
\filldraw (10,0) circle (2.5pt) node[below,yshift=-0.2pt] {\footnotesize C};

\draw [->] (2.05,0.2) to [bend left=40] node[below] {} (2.45,0.2);
\draw [->] (4.55,0.2) to [bend left=40] node[below] {} (4.95,0.2);
\draw [->] (7.05,0.2) to [bend left=40] node[below] {} (7.45,0.2);
\draw [->] (9.55,0.2) to [bend left=40] node[below] {} (9.95,0.2);
\end{scope}
\end{tikzpicture}
\end{center}
\caption{Illustration of F-relaxation (top) and C-relaxation (bottom).}
\label{fig:fc_relaxation}
\end{figure}
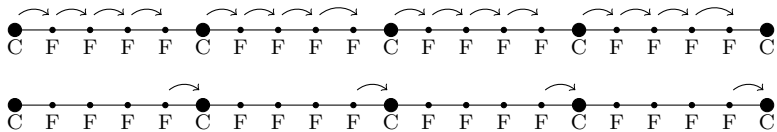

As the parareal method, the two-level MGRIT algorithm can be seen as a multigrid reduction (MGR) method that combines a coarse time stepping method with (parallel) fine time stepping within each coarse time interval. Here, the time stepping from a coarse point $C$ to all neighbouring fine points is also referred to as $F$-relaxation \cite{Falgout2014}. On the other hand, time stepping to a $C$ point from the previous $F$ point is referred to as $C$-relaxation. It should be noted that both types of relaxation are highly parellel and can be combined leading to so-called $CF$- or $FCF$-relaxation. Figure \ref{fig:fc_relaxation} illustrates both relaxation methods. 

\subsection*{Multilevel MGRIT method}

Next, we consider the true multilevel MGRIT method. First, we define a hierarchy of $L$ time discretization meshes, where the time step size for the discretization at level $l \ (l = 0,1, \ldots L)$ is given by $\Delta t_F m^l$. The total number of levels $L$ is related to the coarsening factor $m$ and the total number of fine steps $\Delta t_F$ by $L = log_m(N_t)$. 
Let $\mathbf{A}^{(l)} \mathbf{u}^{(l)} = \mathbf{g}^{(l)}$ denote the linear system of equations based on the considered  time step size at level $l$, where $l=0,1,\dots,L$. The MGRIT method can then be written as follows:
b
\begin{tcolorbox}[colback=black!5!white,colframe=black!80!black]
\textbf{MGRIT} \\
\textbf{if} l == L 
\begin{itemize}
	\item Solve $\mathbf{A}^{(L)} \mathbf{u}^{(L)} = \mathbf{g}^{(L)}$
\end{itemize}
\textbf{else}
\begin{itemize}
	\item Apply FCF-relaxation on $\mathbf{A}^{(l)} \mathbf{u}^{(l)} = \mathbf{g}^{(l)}$
	\item Restrict the residual $\mathbf{g}^{(L)}-\mathbf{A}^{(l)} \mathbf{u}^{(l)}$ using injection 
	\item Call \textbf{MGRIT} setting $l \rightarrow  l+1$
	\item Update $\mathbf{u}^{(l)} \rightarrow \mathbf{u}^{(l)} + P \mathbf{u}^{(l+1)}$
\end{itemize}
\textbf{end}
\end{tcolorbox}

Here, the prolongation operator $P$ is based on ordering the $F$-points and $C$-points, starting with the $F$-points. The matrix $\mathbf{A}$ can then be written as follows:
\begin{eqnarray}
\mathbf{A} =  \begin{bmatrix}  \mathbf{A}_{FF} &  \mathbf{A}_{FC} \\ \mathbf{A}_{CF} & \mathbf{A}_{CC}\end{bmatrix}. 
\end{eqnarray} 

and the operator $P$ is then defined as the ``ideal interpolation" \cite{Falgout2014}:
\begin{eqnarray}
P = \begin{bmatrix} - \mathbf{A}_{FF} \mathbf{A}_{FC} \\ \mathbf{I}_C \end{bmatrix}.
\end{eqnarray}
The recursive algorithm described above leads to a so-called $V$-cycle. However, as with standard multigrid methods, alternative cycle types (i.e. $W$-cycles, $F$-cycles) can be defined. At all levels of the multigrid hierarchy, the operators are obtained by rediscretizing Equation \eqref{eq:heat} using a different time step size. 

\section{$p$-multigrid method}
\label{sec:pmg}

Within the MGRIT algorithm, the action of $\mathbf{A}_{\Delta}$ is computed in parallel to iteratively solve the coarse system as described in Equation \eqref{eq:mgrit_coarse}. Assuming a Backward-Euler time integration scheme, the following linear system of equations is solved within each time interval at every iteration:
\begin{eqnarray} 
\left ( \mathbf{M} + \kappa \Delta t \mathbf{K} \right ) \mathbf{u}^{k+1}= \mathbf{M} \mathbf{u}^{k} + \Delta t \mathbf{f}, \ k = 0,\ldots,N_t. 
\end{eqnarray}

In a recent paper by the authors \cite{Tielen2021}, this linear system of equations was solved within MGRIT by means of a (diagonally preconditioned) Conjugate-Gradient method. However, as the condition number of the system matrix increases exponentially in IgA with the approximation order $p$, the use of standard iterative solvers becomes less efficient for higher values of $p$. As a consequence, alternative solution techniques have been developed in recent years to overcome this dependency \cite{Donatelli2017,Hofreither2017,Riva2019,Sogn2019}. 

In this paper, a $p$-multigrid method \cite{Tielen2019} using an ILUT smoother will be adopted to solve the linear systems within MGRIT. Within the $p$-multigrid method, a low-order correction is obtained (at level $p=1$) to update the solution at the high-order level. Starting from the high-order problem, the following steps are performed \cite{Tielen2019}:

\begin{enumerate}
\item Apply a f\/ixed number $\nu_1$ of presmoothing steps to the initial guess $\mathbf{u}_{h,p}^{(0,0)}$:
\begin{eqnarray} \label{smooth}
\mathbf{u}_{h,p}^{(0,m)} = \mathbf{u}_{h,p}^{(0,m-1)} + \mathcal{S}_{h,p} \left (\mathbf{f}_{h,p} - \mathbf{A}_{h,p} \mathbf{u}_{h,p}^{(0,m-1)} \right ),  \hspace{0.25cm} m=1,\ldots, \nu_1, 
\end{eqnarray}
where $\mathcal{S}_{h,p}$ is a smoothing operator applied to the high-order problem. 
\item Determine the residual at level $p$ and project it onto the space $\mathcal{V}_{h,1}$ using the restriction operator $\mathcal{I}_{p}^{1}$:
\begin{eqnarray}
\mathbf{r}_{h,1} = \mathcal{I}_{p}^{1} \left (  \mathbf{f}_{h,p} - \mathbf{A}_{h,p} \mathbf{u}_{h,p}^{(0,\nu_1)} \right ).
\end{eqnarray} 
\item Solve the residual equation to determine the coarse grid error: 
\begin{eqnarray} \label{eq:res}
\mathbf{A}_{h,1} \mathbf{e}_{h,1} = \mathbf{r}_{h,1}. 
\end{eqnarray}
\item Project the error $\mathbf{e}_{h,1}$ onto the space $\mathcal{V}_{h,p}$ using the prolongation operator $\mathcal{I}_{1}^p$ and update $\mathbf{u}_{h,p}^{(0,\nu_1)}$:
\begin{eqnarray}
\mathbf{u}_{h,p}^{(0,\nu_1)} := \mathbf{u}_{h,p}^{(0,\nu_1)} + \mathcal{I}_{1}^p \left (\mathbf{e}_{h,1} \right ).
\end{eqnarray}
\item Apply $\nu_2$ postsmoothing steps of the form~\eqref{smooth} on the updated solution to obtain $\mathbf{u}_{h,p}^{(0,\nu_1 + \nu_2)} =: \mathbf{u}_{h,p}^{(1,0)}$.
\end{enumerate}

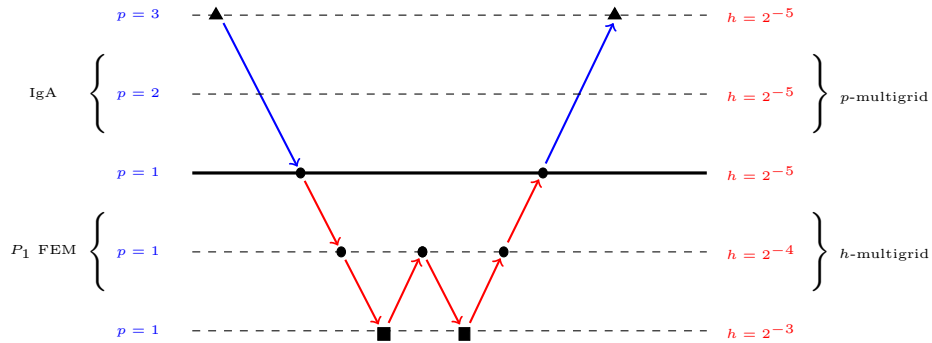
\begin{figure}[h!]
\begin{center}
\begin{tikzpicture}[xscale = 0.36, yscale = 0.42]
\draw (-14,4.5) node {\tiny $\textcolor{blue}{p=3}$};
\draw (9,4.5) node {\tiny $\textcolor{red}{h = 2^{-5}}$};
\draw (-14, 2) node {\tiny $\textcolor{blue}{p=2}$};
\draw (9, 2) node {\tiny $\textcolor{red}{h = 2^{-5}}$};  
\draw (-14,-0.5) node {\tiny $\textcolor{blue}{p=1}$};
\draw (9, -0.5) node {\tiny $\textcolor{red}{h = 2^{-5}}$};
\draw (-14,-3.0) node {\tiny $\textcolor{blue}{p=1}$};
\draw (9, -3.0) node {\tiny $\textcolor{red}{h = 2^{-4}}$};
\draw (-14,-5.5) node {\tiny $\textcolor{blue}{p=1}$};
\draw (9, -5.5) node {\tiny $\textcolor{red}{h = 2^{-3}}$};
\draw (13, 2) node {\Bigg \} \tiny $p$-multigrid};
\draw (13, -3) node {\Bigg \} \tiny $h$-multigrid};
\draw (-15.5, 2) node { \Bigg\{ };
\draw (-15.5, -3) node {\Bigg \{};
\draw (-17.5, 2) node {\tiny IgA};
\draw (-17.5, -3) node {\tiny $P_1$ FEM};

\draw [->,thick,blue] (-10.85,4.25) --  node[sloped, anchor=center, below]{}(-8.12,-0.3);
\draw [->,thick,red] (-7.85,-0.75) --  node[sloped, anchor=center, below]{}(-6.62,-2.8);
\draw [->,thick,red] (-6.35,-3.25) -- node[sloped, anchor=center, below]{}(-5.12,-5.3);
\draw [->,thick,red] (-4.75,-5.25) -- node[sloped, anchor=center, below]{} (-3.62,-3.2);

\draw [->,thick,red] (-3.35,-3.25) -- node[sloped, anchor=center, below]{}(-2.12,-5.3);
\draw [->,thick,red] (-1.75,-5.25) -- node[sloped, anchor=center, below]{} (-0.62,-3.2);
\draw [->,thick,red] (-0.35,-2.7) --  node[sloped, anchor=center, below]{}(0.80,-0.7);
\draw [->,thick,blue] (1.05,-0.2)  --  node[sloped, anchor=center, below]{}(3.55,4.25);

\draw[dashed] (-12, 4.5) -- (7, 4.5);
\draw[dashed] (-12, 2) -- (7, 2);
\draw[very thick] (-12, -0.5) -- (7, -0.5);
\draw[dashed] (-12, -3.0) -- (7, -3.0);
\draw[dashed] (-12, -5.5) -- (7, -5.5);

\node[scale=1.35] at (-11.125,4.5) {\pgfuseplotmark{triangle*}};
\filldraw (-8,-0.5) circle (4.5pt) node[]{};
\draw[black, fill] (-5.15,-5.8) rectangle (-4.70,-5.4);
\filldraw (-6.5,-3.0) circle (4.5pt) node[]{};
\filldraw (-3.5,-3.0) circle (4.5pt) node[]{};
\draw[black, fill] (-2.15,-5.8) rectangle (-1.75,-5.4);
\filldraw (-3.5,-3.0) circle (4.5pt) node[]{};
\filldraw (-0.5,-3.0) circle (4.5pt) node[]{};
\filldraw (0.95,-0.5) circle (4.5pt) node[]{};
\node[scale=1.35] at (3.60,4.5) {\pgfuseplotmark{triangle*}};
\end{tikzpicture}
\end{center}
\caption{Illustration of the $p$-multigrid method \cite{Tielen2019}. At $p=1$, Gauss-Seidel is adopted as a smoother ($\bullet$), whereas at the high-order level ILUT is applied ($\blacktriangle$). At the coarsest level, a direct solver is applied to solve the residual equation ($\blacksquare$).}
\label{fig:pmg}
\end{figure}

To approximately solve the residual equation given by Equation \eqref{eq:res} a single W-cycle of a standard $h$-multigrid method \cite{Hackbush1985}, using canonical prolongation and weighted restriction, is applied. As the level $p=1$ corresponds to a low-order Lagrange discretization, an $h$-multigrid method (using Gauss-Seidel as a smoother) is known to be both efficient and cheap \cite{Hackbush1985,Oosterlee2001}. The resulting $p$-multigrid adopted throughout this paper is shown in Figure \ref{fig:pmg}.

Note that, we directly restrict the residual at the high-order level to level $p=1$. This aggresive $p$-coarsening strategy has shown to significantly improve the computational efficiency of the resulting $p$-multigrid method \cite{Tielen2021b}.

Prolongation and restriction operators based on an $L_2$ projection are adopted to transfer vectors from the high-order level to the low-order level (and vice versa). The operators have been used extensively in the literature \cite{Briggs2000,Brenner1994,Sampath2010} and are given by:
\begin{eqnarray} \label{prolongation}
\mathcal{I}_{1}^p(\mathbf{{}v}_{1}) = (\mathbf{M}_p)^{-1} \mathbf{P}_{1}^{p} \ \mathbf{v}_{1}, \qquad 
\mathcal{I}_{p}^{1}(\mathbf{v}_{p}) = (\mathbf{M}_{1})^{-1} \mathbf{P}_{p}^{{1}} \ \mathbf{v}_{p}.
\end{eqnarray}
Here, the mass matrix $\mathbf{M}_p$ and transfer matrix $\mathbf{P}_{1}^{p}$ are defined as follows:
\begin{eqnarray} \label{PM}
 ( \mathbf{M}_{p})_{(i,j)} := \int_{\Omega} \Phi_{i,p} \Phi_{j,p}  \hspace{0.1cm} \text{d} \Omega,  \hspace{1.0cm} (\mathbf{P}^{p}_{1})_{(i,j)} :=  \int_{\Omega} \Phi_{i,p} \Phi_{j,1}  \hspace{0.1cm} \text{d} \Omega.
\end{eqnarray}

To prevent the explicit solution of a linear system of equations for each projection step, the consistent mass matrix $\mathbf{M}$ in both transfer operators is replaced by its lumped counterpart $\mathbf{M}^L$ by applying row-sum lumping (i.e. $\mathbf{M}_{(i,i)}^L = \sum_{j=1}^{N_{\rm dof}} \mathbf{M}_{(i,j)}$). Note that, row-sum lumping can be applied within the variational formulation, due to the partition of unity and non-negativity of the B-spline basis functions. \\
\\
Various choices can be made with respect to the smoother. The use of Gauss-Seidel or (damped) Jacobi as a smoother at level $p$ leads to convergence rates that depend significantly on the approximation order $p$ \cite{Tielen2017}. Alternative smoothers have been developed in recent years to overcome this shortcoming \cite{Donatelli2017,Hofreither2017,Riva2019,Sogn2019}. In particular, the use of ILUT factorizations has shown to be very effective and will therefore be adopted throughout the remainder of this paper.

\section{Numerical results}
\label{sec:results}

\begin{figure}[h]
\centering
\begin{subfigure}{.4\linewidth}
\centering
\includegraphics[scale=0.21]{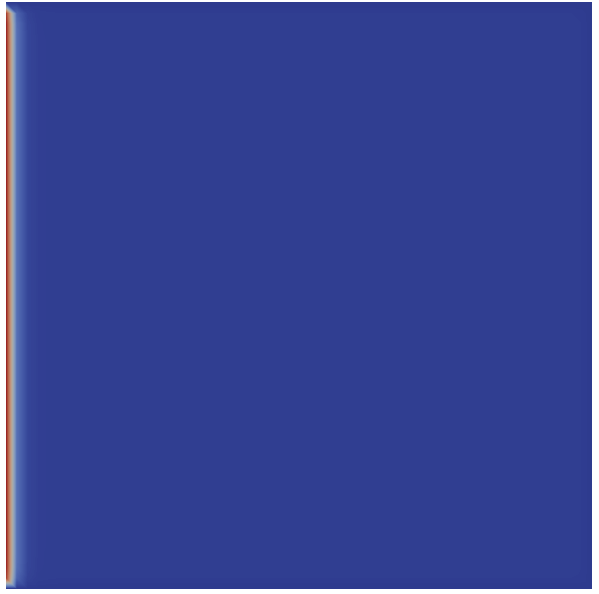}
\caption{$T=0$}
\end{subfigure}
\begin{subfigure}{.4\linewidth}
\centering
\includegraphics[scale=0.21]{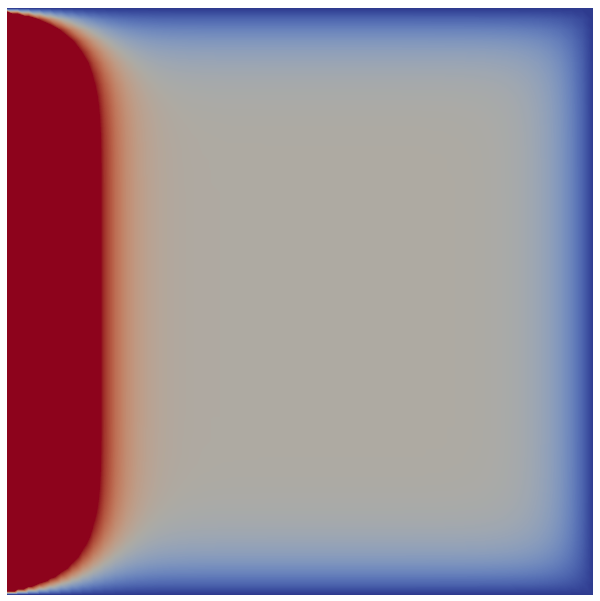}
\caption{$T=0.005$}
\end{subfigure}
\begin{subfigure}{.4\linewidth}
\centering
\includegraphics[scale=0.21]{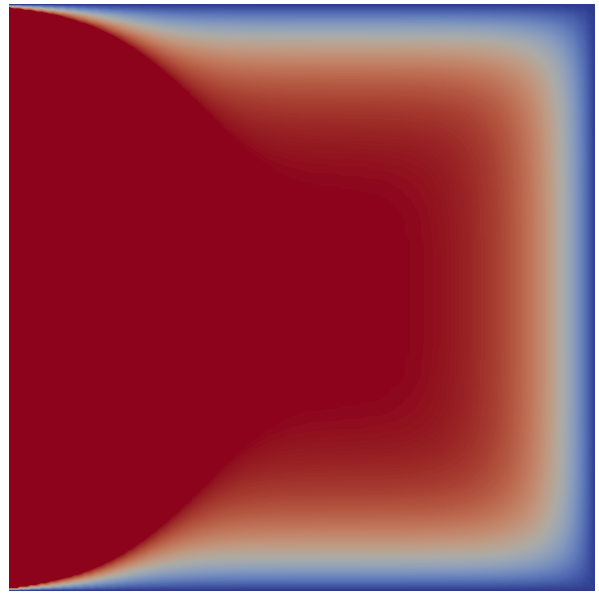}
\caption{$T=0.010$}
\end{subfigure}
\begin{subfigure}{.4\linewidth}
\centering
\includegraphics[scale=0.21]{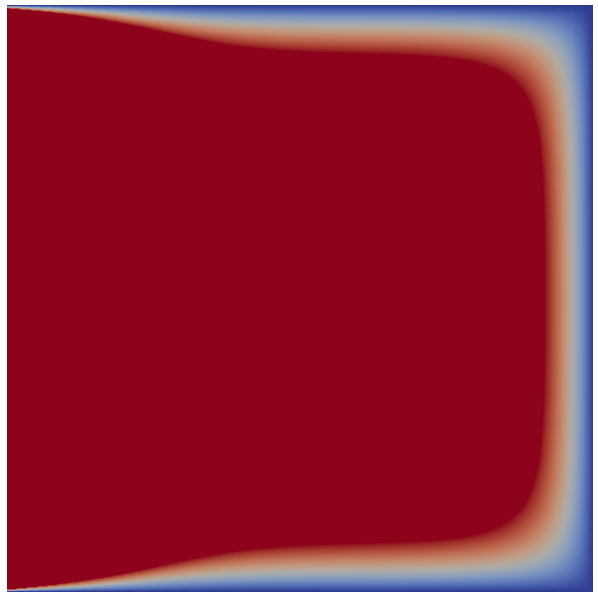}
\caption{$T=0.020$}
\end{subfigure}
\caption{Solution to MP-1 at different times $T$ using a inhomogeneous Neumann boundary condition at the left boundary using quadratic B-spline basis functions.}
\label{fig:solution}
\end{figure}

To assess the quality of MGRIT when applied in combination with a $p$-multigrid method within Isogeometric Analysis, we consider the time-dependent heat equation in two dimensions given by Equation \eqref{eq:heat}. Figure \ref{fig:solution} shows the resulting solution $u$ at different time instances for $\Omega = [0,1]^2$. Here, an inhomogeneous Neumann boundary condition is applied at the left boundary. 

Based on a spatial discretization with B-spline basis functions of order $p$ and mesh width $h$, MGRIT is applied to iteratively solve the resulting equation. In particular, we will investigate the parallel performance of MGRIT on modern computer architectures. The open-source C++ library G+Smo \cite{gismo} is used to discretize the model problem in space, while, for the MGRIT algorithm, the parallel-in-time code XBraid, developed at Lawrence Livermore National Lab, is adopted \cite{xbraid}.

\begin{figure}[h!]
\centering
\begin{tikzpicture}[scale=0.7]
\draw (0,0) -- (4,0) -- (4,4) -- (0,4) -- (0,0);

\draw (8,0) -- (10,0);
\draw (6,4) -- (6,2);
\draw (10,0) arc (0:90:4cm);
\draw (8,0) arc (0:90:2cm);

\draw[dashed] (14,0) -- (14,4);
\draw[dashed] (12,2) -- (16,2);
\draw (12,0) -- (16,0) -- (16,4) -- (12,4) -- (12,0);
\end{tikzpicture}
\caption{Spatial domains $\Omega$ considered in \cite{Tielen2021}.}
\label{fig:domain}
\end{figure}
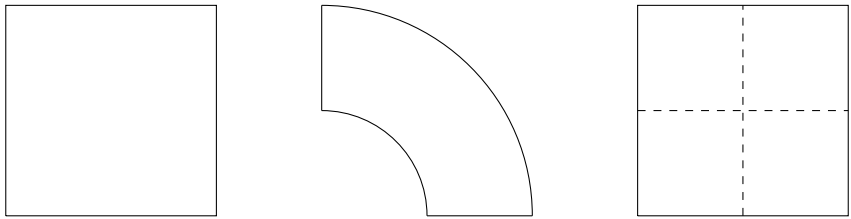

As a model problem, we solve Equation \eqref{eq:heat} on the time domain $T=[0,0.1]$ with $\kappa =1$. Furthermore, the right-hand side is chosen equal to one and the initial condition is equal to zero. The MGRIT method is said to have reached convergence if the relative residual (in the $L_2$ norm) at the end of an iteration is smaller or equal to $10^{-10}$, unless stated otherwise.
\newpage
As a starting point, we briefly summarize the results obtained in a previous paper of the authors (see \cite{Tielen2021}). There, numerical results were obtained for same model problem using different hierarchies (i.e. a V-cycle, F-cycle and two-level method), time integration schemes (i.e. backward Euler, forward Euler and Crank-Nicolson) and domains of interest (see Figure \ref{fig:domain}). 

In general, it was observed that MGRIT converged in a low number (i.e. $5 \sim 10$) of iterations, although the number of iterations was slightly higher when V-cycles were adopted instead of F-cycles or a two-level method. Furthermore, the number of iterations was independent of the mesh width $h$, approximation order of the B-spline basis functions $p$ and the number of time steps $N_t$ for all considered hierarchies and domains of interest. As expected from sequential time stepping methods, the use of the implicit backward Euler within MGRIT lead to the most stable time integration method. Finally, CPU times were obtained for a limited number of processors, showing a strong dependency on the approximation order $p$ when the Conjugate Gradient method was applied as a spatial solver within MGRIT. 

In this section, we investigate the effect of using a $p$-multigrid method for the spatial solves compared to the use of a Conjugate Gradient method. Furthermore, we investigate the weak and strong scaling of MGRIT on modern architectures when applied in the context of IgA.

\subsection*{Iteration numbers}
As a first step, we compare the number of MGRIT iterations to reach convergence when a $p$-multigrid method is adopted while keeping all other parameters the same. Table \ref{tab:7} shows the results when the mesh width is kept constant ($h=2^{-6}$) for the unit square and a quarter annulus as our domain $\Omega$ when adopting V-cycles. For both benchmarks and all configurations, the number of iterations needed with MGRIT to reach convergence is identical, which was observed as well in \cite{Tielen2021} in case a Conjugate Gradient method was used for the spatial solves. 

\begin{table}[h]
\centering
\begin{tabular}{l|cccc|cccc|}
  \multicolumn{1}{c|}{}   & \multicolumn{4}{c|}{Unit Square}      & \multicolumn{4}{c|}{Quarter Annulus}    \\  
                         & $p=2$     & $p=3$ & $p=4$ & $p=5$ & $p=2$  & $p=3$    & $p=4$  & $p=5$  \\ \hline
        $N_t=250$        &  $10$     & $10$  & $10$  &  $10$  &  $10$  & $10$    &  $10$  &   $10$     \\
        $N_t=500$        &  $10$     & $10$  & $10$  &  $10$  &  $10$  & $10$    &  $10$  &   $10$     \\
        $N_t=1000$       &  $11$     & $11$  & $11$  &  $11$  &  $11$  & $11$    &  $11$  &   $11$     \\
        $N_t=2000$       &  $11$     & $11$  & $11$  &  $11$  &  $11$  & $11$    &  $11$  &   $11$     \\
\end{tabular}
\caption{Number of MGRIT iterations for solving Equation \eqref{eq:heat} on the unit square and a quarter annulus when adopting V-cycles for a varying number of time steps. Here $p$-multigrid is adopted for the spatial solves.}
\label{tab:7}
\end{table}

Table \ref{tab:8} shows the results when the number of time steps is kept constant ($N_t=100$) for both benchmarks when adopting V-cycles. Results show a similar number of iterations compared to the use of a Conjugate Gradient method.

\begin{table}[h]
\centering
\begin{tabular}{l|cccc|cccc|}
  \multicolumn{1}{c|}{}   & \multicolumn{4}{c|}{Unit Square}      & \multicolumn{4}{c|}{Quarter Annulus}    \\
                         & $p=2$ & $p=3$ & $p=4$ & $p=5$  & $p=2$  & $p=3$   & $p=4$  & $p=5$  \\ \hline
        $h=2^{-6}$       &  $9$  & $9$   & $9$   &  $9$  &  $9$  & $9$   &  $9$  &   $9$     \\
        $h=2^{-7}$       &  $9$  & $9$   & $9$   &  $9$  &  $9$  & $9$   &  $9$  &   $9$     \\
        $h=2^{-8}$       &  $10$  & $10$   & $10$   &  $10$  &  $9$   & $9$    &  $9$   &   $9$     \\
        $h=2^{-9}$       &  $10$  & $10$   & $10$   &  $10$  &  $10$  & $10$   &  $10$  &   $10$     \\
\end{tabular}
\caption{Number of MGRIT iterations for solving Equation \eqref{eq:heat} on the unit square and a quarter annulus when adopting V-cycles for varying mesh widths. Here $p$-multigrid is adopted for the spatial solves.}
\label{tab:8}
\end{table}

\subsection*{CPU timings}

CPU timings have been obtained when a $p$-multigrid method or Conjugate Gradient method is adopted for the spatial solves within MGRIT. As in the previous section, we adopt V-cycles, a mesh width of $h=2^{-6}$ and the unit square as our domain of interest. Note that the corresponding iteration numbers can be found in Table \ref{tab:7}. The computations are performed on three nodes, which consist each of an Intel(R) i7-10700 (@ 2.90GHz) processor with $8$ cores.

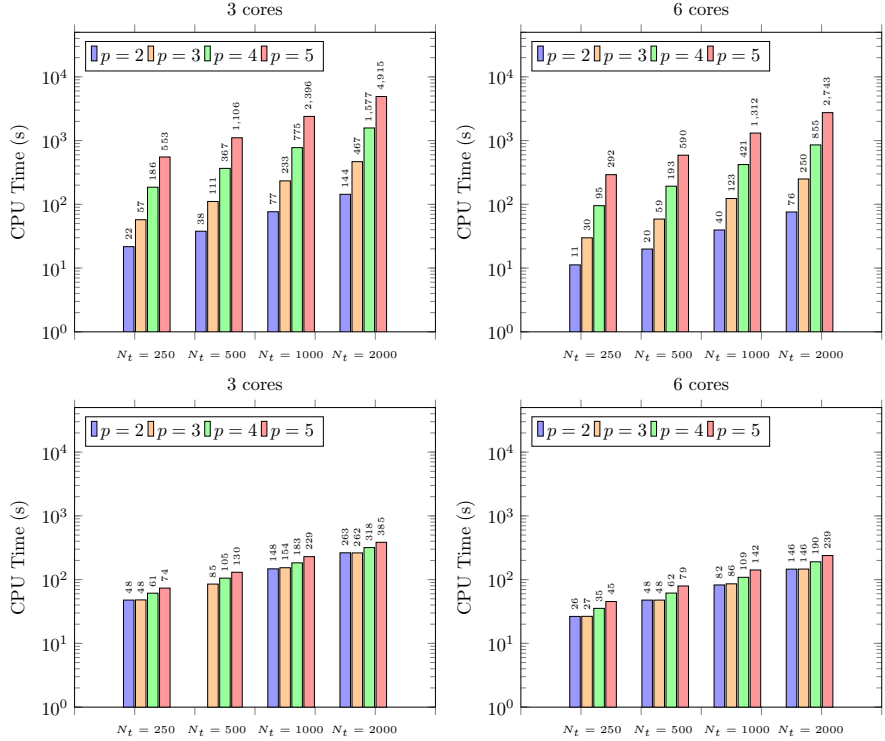
\begin{figure}[h!]
\centering
\begin{tikzpicture}[scale=0.7]
    \begin{axis}[nodes near coords={\pgfmathprintnumber[fixed,fixed zerofill,precision=0]{\pgfplotspointmeta}},point meta=rawy,ybar,every node near coord/.append style={rotate=90, anchor=west,font=\tiny},legend columns=4,legend pos=north west, xmin= 0,xmax=30,ymode=log,ymin=1,ymax=50000,ylabel=CPU Time (s),line width=0.25pt,bar width=0.2cm,bar shift=8mm,xticklabels=\emptynodes,title=$3$ cores]
       \addplot[fill=blue!40] plot coordinates 
      {
        (1,    21.68)
        (7,    37.87)
        (13,   76.53)
        (19,   143.73)
       };
      \addlegendentry{$p=2$}
       \addplot[fill=orange!40] plot coordinates 
      {
        (2,    57.34)
        (8,    110.89)
        (14,   233.34)
        (20,   467.00)        
       };
      \addlegendentry{$p=3$}
       \addplot[fill=green!40] plot coordinates 
      {
        (3,    185.65)
        (9,    366.51)
        (15,   774.97)
        (21,   1576.87)
       };
      \addlegendentry{$p=4$}
      \addplot[fill=red!40] plot coordinates 
      {
        (4,    553.23)
        (10,   1106.46)
        (16,   2395.64)
        (22,   4915.37)
        };
      \addlegendentry{$p=5$}
      \end{axis} 
        \draw (1.35,-0.45)  node{\scalebox{.8}{\tiny $N_t=250$}};
        \draw (2.70,-0.45)  node{\scalebox{.8}{\tiny $N_t=500$}};
        \draw (4.10,-0.45)  node{\scalebox{.8}{\tiny $N_t=1000$}};
        \draw (5.50,-0.45)  node{\scalebox{.8}{\tiny $N_t=2000$}};
\end{tikzpicture}
  \begin{tikzpicture}[scale=0.7]
    \begin{axis}[nodes near coords={\pgfmathprintnumber[fixed,fixed zerofill,precision=0]{\pgfplotspointmeta}},point meta=rawy,ybar,every node near coord/.append style={rotate=90, anchor=west,font=\tiny},legend columns=4,legend pos=north west, xmin= 0,xmax=30,ymode=log,ymin=1,ymax=50000,ylabel=CPU Time (s),line width=0.25pt,bar width=0.2cm,bar shift=8mm,xticklabels=\emptynodes,title=$6$ cores]
       \addplot[fill=blue!40] plot coordinates 
      {
        (1,    11.23)
        (7,    19.90)
        (13,   39.59)
        (19,   75.87)
       };
      \addlegendentry{$p=2$}
       \addplot[fill=orange!40] plot coordinates 
      {
        (2,    29.84)
        (8,    58.76)
        (14,   123.37)
        (20,   249.75)        
       };
      \addlegendentry{$p=3$}
       \addplot[fill=green!40] plot coordinates 
      {
        (3,    95.30)
        (9,    192.63)
        (15,   420.53)
        (21,   855.00)
       };
      \addlegendentry{$p=4$}
      \addplot[fill=red!40] plot coordinates 
      {
        (4,    292.04)
        (10,   590.29)
        (16,   1311.94)
        (22,   2743.36)
        };
      \addlegendentry{$p=5$}
      \end{axis}
       \draw (1.35,-0.45)  node{\scalebox{.8}{\tiny $N_t=250$}};
        \draw (2.70,-0.45)  node{\scalebox{.8}{\tiny $N_t=500$}};
        \draw (4.10,-0.45)  node{\scalebox{.8}{\tiny $N_t=1000$}};
        \draw (5.50,-0.45)  node{\scalebox{.8}{\tiny $N_t=2000$}};
\end{tikzpicture}
\begin{tikzpicture}[scale=0.7]
    \begin{axis}[nodes near coords={\pgfmathprintnumber[fixed,fixed zerofill,precision=0]{\pgfplotspointmeta}},point meta=rawy,ybar,every node near coord/.append style={rotate=90, anchor=west,font=\tiny},legend columns=4,legend pos=north west, xmin= 0,xmax=30,ymode=log,ymin=1,ymax=50000,ylabel=CPU Time (s),line width=0.25pt,bar width=0.2cm,bar shift=8mm,xticklabels=\emptynodes,title=$3$ cores]
       \addplot[fill=blue!40] plot coordinates 
      {
        (1,    47.76)
        (13,   147.64)
        (19,   263.36)
       };
      \addlegendentry{$p=2$}
       \addplot[fill=orange!40] plot coordinates 
      {
        (2,    48.06)
        (8,    85.13)
        (14,   153.60)
        (20,   262.39)        
       };
      \addlegendentry{$p=3$}
       \addplot[fill=green!40] plot coordinates 
      {
        (3,    61.27)
        (9,    105.47)
        (15,   183.16)
        (21,   317.94)
       };
      \addlegendentry{$p=4$}
      \addplot[fill=red!40] plot coordinates 
      {
        (4,    73.60)
        (10,   130.42)
        (16,   229.21)
        (22,   384.98)
        };
      \addlegendentry{$p=5$}
      \end{axis}  
        \draw (1.35,-0.45)  node{\scalebox{.8}{\tiny $N_t=250$}};
        \draw (2.70,-0.45)  node{\scalebox{.8}{\tiny $N_t=500$}};
        \draw (4.10,-0.45)  node{\scalebox{.8}{\tiny $N_t=1000$}};
        \draw (5.50,-0.45)  node{\scalebox{.8}{\tiny $N_t=2000$}};
\end{tikzpicture}
  \begin{tikzpicture}[scale=0.7]
    \begin{axis}[nodes near coords={\pgfmathprintnumber[fixed,fixed zerofill,precision=0]{\pgfplotspointmeta}},point meta=rawy,ybar,every node near coord/.append style={rotate=90, anchor=west,font=\tiny},legend columns=4,legend pos=north west, xmin= 0,xmax=30,ymode=log,ymin=1,ymax=50000,ylabel=CPU Time (s),line width=0.25pt,bar width=0.2cm,bar shift=8mm,xticklabels=\emptynodes,title=$6$ cores]
       \addplot[fill=blue!40] plot coordinates 
      {
        (1,    26.48)
        (7,    47.75)
        (13,   82.48)
        (19,   146.02)
       };
      \addlegendentry{$p=2$}
       \addplot[fill=orange!40] plot coordinates 
      {
        (2,    26.55)
        (8,    47.73)
        (14,   85.88)
        (20,   146.40)        
       };
      \addlegendentry{$p=3$}
       \addplot[fill=green!40] plot coordinates 
      {
        (3,    35.43)
        (9,    61.52)
        (15,   108.69)
        (21,   190.27)
       };
      \addlegendentry{$p=4$}
      \addplot[fill=red!40] plot coordinates 
      {
        (4,    45.36)
        (10,   79.49)
        (16,   141.60)
        (22,   238.76)
        };
      \addlegendentry{$p=5$}
      \end{axis}
        \draw (1.35,-0.45)  node{\scalebox{.8}{\tiny $N_t=250$}};
        \draw (2.70,-0.45)  node{\scalebox{.8}{\tiny $N_t=500$}};
        \draw (4.10,-0.45)  node{\scalebox{.8}{\tiny $N_t=1000$}};
        \draw (5.50,-0.45)  node{\scalebox{.8}{\tiny $N_t=2000$}};
\end{tikzpicture}
\caption{CPU times for MGRIT using V-cycles and backward Euler on the unit square for a fixed problem size ($h=2^{-6}$) adopting a different number of processors. Here the Conjugate Gradient method (top) and a $p$-multigrid method (bottom) are used for all spatial solves within MGRIT.}
\label{fig:cpu1}
\end{figure}

Figure \ref{fig:cpu1} shows the CPU time needed to reach convergence for a varying number of cores, a different number of time steps and different values of $p$. When the Conjugate Gradient method is adopted for the spatial solves, doubling the number of time steps leads to an increase of the CPU time by a factor of two. Furthermore, it can be observed that the CPU times significantly increase for higher values of $p$ which is related to the spatial solves required at every time step. As standard iterative solvers (like the Conjugate Gradient method) have a detoriating performance for increasing values of $p$, more iterations are required to reach convergence for each spatial solve, resulting in higher computational costs of the MGRIT method. When focussing on the number of cores, it can be seen that doubling the number of cores significantly reduces the CPU time needed to reach convergence. More precisely, a reduction of $45-50\%$ can be observed when doubling the number of cores to $6$, implying the MGRIT algorithm is highly parallelizable. 

As with the use of the Conjugate Gradient method, doubling the number of time steps leads to an increase of the CPU time by a factor of two when a $p$-multigrid method is adopted. However, the dependency of the CPU times on the approximation order is significantly mitigated, which leads to a serious decrease of the CPU times compared to the use of the Conjugate Gradient method when higher values of $p$ are considered. Again, increasing the number of cores from $3$ to $6$, reduces the CPU time needed to reach convergence by $45-50\%$. 

These results indicate that MGRIT combined with a $p$-multigrid method leads to an overall more efficient method. Therefore, a large computer cluster will be considered in the next section to further investigate the scalability of MGRIT (i.e. weak and strong scalability) when combined with a $p$-multigrid method within IgA.

\section{Scalability}
\label{sec:scalable}
In the previous sections, we applied MGRIT adopting a relatively low number of cores. Here, it is shown that the use of a $p$-multigrid method significantly reduces the dependency of the CPU timings on the approximation order. In this section, we investigate the scalability of MGRIT (combined with a $p$-multigrid method) on a modern architecture. More precisely, we will investigate both strong and weak scalability on the Lisa system, one of the nationally used clusters of the Netherlands \footnote{https://userinfo.surfsara.nl/systems/lisa}. 

\subsection*{Strong scalability}
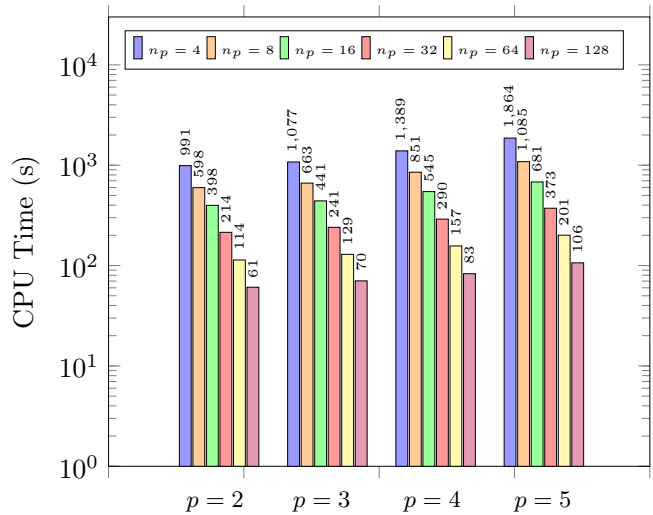
\begin{figure}[h!]
\centering
\begin{tikzpicture}[xscale=1.05,yscale=1.05]
    \begin{axis}[nodes near coords={\pgfmathprintnumber[fixed,fixed zerofill,precision=0]{\pgfplotspointmeta}},point meta=rawy,ybar,every node near coord/.append style={rotate=90, anchor=west,font=\tiny},legend columns=6,legend pos=north west, xmin= 0,xmax=40,ymode=log,ymin=1,ymax=30000,ylabel=CPU Time (s),line width=0.25pt,bar width=0.15cm,bar shift=8mm,xticklabels=\emptynodes]
       \addplot[fill=blue!40] plot coordinates 
      {
        (1,    990.84)
        (9,    1076.99)
        (17,   1389.17)
        (25,   1864.40)
       };
      \addlegendentry{\scalebox{0.9}{\tiny $n_p = 4$}}
       \addplot[fill=orange!40] plot coordinates 
      {
        (2,    597.50)
        (10,    663.02)
        (18,   851.42)
        (26,   1084.54)        
       };
      \addlegendentry{\scalebox{0.9}{\tiny $n_p = 8$}}
       \addplot[fill=green!40] plot coordinates 
      {
        (3,    398.24)
        (11,    441.07)
        (19,   545.46)
        (27,   680.64)
       };
      \addlegendentry{\scalebox{0.9}{\tiny $n_p = 16$}}
      \addplot[fill=red!40] plot coordinates 
      {
        (4,    214.41)
        (12,   240.52)
        (20,   290.21)
        (28,   372.74)
        };
      \addlegendentry{\scalebox{0.9}{\tiny $n_p = 32$}}
      \addplot[fill=yellow!40] plot coordinates 
      {
        (5,    113.73)
        (13,   129.36)
        (21,   156.78)
        (29,   200.77)
        };
      \addlegendentry{\scalebox{0.9}{\tiny $n_p = 64$}}     
      \addplot[fill=purple!40] plot coordinates 
      {
        (6,    60.84)
        (14,   70.38)
        (22,   82.91)
        (30,   106.42)
        };
      \addlegendentry{\scalebox{0.9}{\tiny $n_p = 128$}}
      \end{axis} 
        \draw (1.35,-0.45)  node{\small $p=2$};
        \draw (2.70,-0.45)  node{\small $p=3$};
        \draw (4.10,-0.45)  node{\small $p=4$};
        \draw (5.50,-0.45)  node{\small $p=5$};
\end{tikzpicture}
\caption{Strong scalability study for MGRIT using V-cycles and backward Euler on the unit square. Here $p$-multigrid is used for all spatial solves within MGRIT.}
\label{fig:strongscaling}
\end{figure}

First, we fix the total problem size and increase the number of cores (i.e. strong scalability). That is, we consider the same benchmark problem as in the previous sections, but with a mesh width of $h=2^{-6}$ and a number of time steps $N_t$ of $10.000$. As before, backward Euler is applied for the time integration and V-cycles are adopted as MGRIT hierarchy. Figure \ref{fig:strongscaling} shows the CPU times needed to reach convergence for a varying number of Intel Xeon Gold 6130 (@ 2.10GHz) processors, where each processor consists of $16$ cores. For all values of $p$, increasing the number of processors $n_p$ leads to significant speed-ups which illustrates the parallizability of the MGRIT method up to $2048$ cores ($128$ processors).        

\newpage

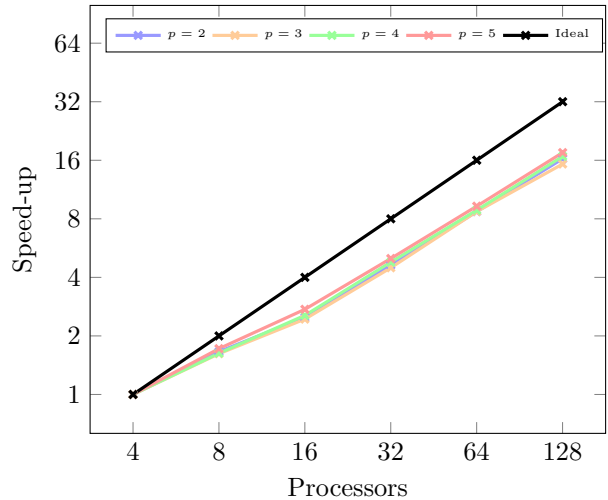
\begin{figure}[h!]
\centering
\begin{tikzpicture}
  \begin{loglogaxis}[xlabel=Processors,ylabel=Speed-up,log basis x=2,xticklabel={\xinttheiexpr[0]2^\tick\relax},log basis y=2,yticklabel={\xinttheiexpr[0]2^\tick\relax},legend columns=5,legend pos=north west,ymax=100]
  \addplot[color=blue!40,mark=x,line width=1.2] coordinates {
    (4, 1)
    (8, 1.66)
    (16,2.49)
    (32,4.62)
    (64,8.71)
    (128,16.29)
  };
  \addlegendentry{\scalebox{0.88}{\tiny $p=2$}}
    \addplot[color=orange!40,mark=x,line width=1.2] coordinates {
    (4, 1)
    (8, 1.62)
    (16,2.44)
    (32,4.48)
    (64,8.71)
    (128,15.30)
  };
  \addlegendentry{\scalebox{0.88}{\tiny $p=3$}}
    \addplot[color=green!40,mark=x,line width=1.2] coordinates {
    (4, 1)
    (8, 1.63)
    (16,2.55)
    (32,4.79)
    (64,8.86)
    (128,16.76)
  };
  \addlegendentry{\scalebox{0.88}{\tiny $p=4$}}
    \addplot[color=red!40,mark=x,line width=1.2] coordinates {
    (4, 1)
    (8, 1.72)
    (16,2.74)
    (32,5.00)
    (64,9.29)
    (128,17.52)
  };
  \addlegendentry{\scalebox{0.88}{\tiny $p=5$}}
  \addplot[color=black,mark=x,line width=1.2] coordinates {
    (4, 1)
    (8, 2)
    (16,4)
    (32,8)
    (64,16)
    (128,32)
  };
  \addlegendentry{\scalebox{0.88}{\tiny Ideal}}
  \end{loglogaxis}
\end{tikzpicture}
\caption{Speed-up with MGRIT using V-cycles and backward Euler on the unit square. Here $p$-multigrid is used for all spatial solves within MGRIT.}
\label{fig:speedup}
\end{figure}

Figure \ref{fig:speedup} shows the obtained speed-ups as a function of the number of processors for different values of $p$ based on the results presented in Figure \ref{fig:strongscaling}. As a comparison, the ideal speed-up has been added, assuming a perfect parallizability of the MGRIT method. The obtained speed-ups remain high, even when the number of processors is further increased to $128$, and is independent of the approximation order $p$.

\subsection*{Weak scalability}
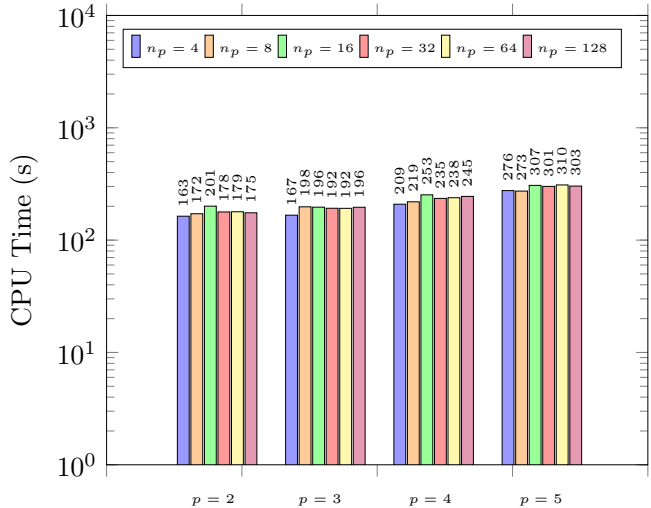
\begin{figure}[h!]
\centering
\begin{tikzpicture}[xscale=1.05,yscale=1.05]
    \begin{axis}[nodes near coords={\pgfmathprintnumber[fixed,fixed zerofill,precision=0]{\pgfplotspointmeta}},point meta=rawy,ybar,every node near coord/.append style={rotate=90, anchor=west,font=\tiny},legend columns=6,legend pos=north west, xmin= 0,xmax=40,ymode=log,ymin=1,ymax=10000,ylabel=CPU Time (s),line width=0.25pt,bar width=0.15cm,bar shift=8mm,xticklabels=\emptynodes]
       \addplot[fill=blue!40] plot coordinates 
      {
        (1,    163.26)
        (9,    166.70)
        (17,   208.75)
        (25,   275.94)
       };
      \addlegendentry{\scalebox{0.9}{\tiny $n_p=4$}}
       \addplot[fill=orange!40] plot coordinates 
      {
        (2,    171.71)
        (10,    197.86)
        (18,   219.28)
        (26,   273.20)        
       };
      \addlegendentry{\scalebox{0.9}{\tiny $n_p=8$}}
       \addplot[fill=green!40] plot coordinates 
      {
        (3,    200.56)
        (11,    196.48)
        (19,   252.95)
        (27,   307.04)
       };
      \addlegendentry{\scalebox{0.9}{\tiny $n_p=16$}}
      \addplot[fill=red!40] plot coordinates 
      {
        (4,    178.16)
        (12,   192.14)
        (20,   235.07)
        (28,   300.81)
        };
      \addlegendentry{\scalebox{0.9}{\tiny $n_p=32$}}
      \addplot[fill=yellow!40] plot coordinates 
      {
        (5,    178.85)
        (13,   191.84)
        (21,   238.45)
        (29,   309.86)
        };
      \addlegendentry{\scalebox{0.9}{\tiny $n_p=64$}}     
      \addplot[fill=purple!40] plot coordinates 
      {
        (6,    175.05)
        (14,   196.13)
        (22,   244.99)
        (30,   302.76)
        };
      \addlegendentry{\scalebox{0.9}{\tiny $n_p=128$}}
      \end{axis} 
        \draw (1.35,-0.45)  node{\tiny $p=2$};
        \draw (2.70,-0.45)  node{\tiny $p=3$};
        \draw (4.10,-0.45)  node{\tiny $p=4$};
        \draw (5.50,-0.45)  node{\tiny $p=5$};
\end{tikzpicture}
\caption{Weak scalability study for MGRIT using V-cycles and backward Euler on the unit square. Here $p$-multigrid is used for all spatial solves within MGRIT.}
\label{fig:weakscaling}
\end{figure}

As a next step, we consider the same benchmark problem but keep the problem size per processor fixed (i.e. weak scalability). In case of four processors, the number of time steps equals $1000$ and is adjusted based on the number of processors. Figure \ref{fig:weakscaling} shows the CPU time needed to reach convergence for a different number of processors and different values of $p$. Clearly, the CPU times remain (more or less) constant when the number of processors is increased, showing the scalability of the MGRIT method. Although the CPU times slightly increase for higher of $p$, the strong $p$-dependency observed with the Conjugate Gradient method is clearly mitigated.  

\newpage

\section{Conclusions}
\label{sec:conclusions}

In this paper, we combined MGRIT with a $p$-multigrid method for discretizations arising in Isogeometric Analysis. Numerical results obtained for two-dimensional benchmark problems show that the use of a $p$-multigrid method for all spatial solves within MGRIT results in CPU times that depend only mildly on the approximation order $p$. This in sharp contrast to standard solvers (e.g. a Conjugate Gradient method), which show a detoriating performance (in terms of CPU times) for higher values of $p$. Furthermore, the obtained CPU times when adopting a $p$-multigrid method are significantly lower for all considered configurations. 

On modern computer architectures, both strong and weak scalability of the resulting MGRIT method have been investigated, showing good scalability up to $2048$ cores. This illustrates the potential of MGRIT (combined with $p$-multigrid) for time-dependent simulations in IgA. Future work will therefore focus on the application of MGRIT to more challening benchmark problems, in particular those where IgA has proven to be a viable alternative to FEM.

\bibliography{mybibfile}

\end{document}